\documentclass[a4paper,12pt,reqno]{amsart}
\usepackage{amsmath,amsfonts,amssymb,amsthm,}
\usepackage{a4,enumerate}
\renewenvironment{proof}{\vskip .2em\noindent{\bfseries Proof.}\ \ }{\hfill$\Box$}

\newtheorem{theorem}{\indent\bf Theorem}
\newtheorem{theoremA}{\indent\bf Theorem}

\newtheorem{lemma}{\indent\bf Lemma}

\theoremstyle{remark}
\newtheorem{remark}{\indent\bf Remark}
\numberwithin{equation}{section}

\title[]{}
\author[]{}
\begin{document}
\begin{center}
{\large\bf ZYGMUND'S TYPE INEQUALITY TO THE POLAR\\[.25em] DERIVATIVE OF A POLYNOMIAL}\\[1em]
{M.S. Pukhta}\\[.35em]
SKUAST-Kashmir, Srinagar 190025, India\\
E-mail: mspukhta\_67@yahoo.co.in
\end{center}

\vskip .25in \thispagestyle{plain} \fontsize{12}{20}\selectfont
\noindent\textbf{Abstract:}
In this paper we improve a result recently proved by Irshad et al.\
[On the Inequalities Concerning to the Polar Derivative of a Polynomial with Restricted Zeroes,
Thai Journal of Mathematics, 2014 (Article in Press)]
and also extend Zygmund's inequality to the polar derivative of a polynomial.



\section{Introduction and Statement of Results}

Let $P(z)$ ba a polynomial of degree $n$, then
\begin{align}
\max\limits_{|z|=1}|P'(z)|\le n\max\limits_{|z|=1} |P(z)|
\end{align}
inequlaity (1.1) is a well known result of S.~Bernstein \cite{1}.
Equality holds in (1.1) if and only if $P(z)$ has all its zeros at the origin.

Inequality (1.1) was extended to $L_p$-norm $p\ge 1$ by Zygmund \cite{2},
who proved that if $P(z)$ is a polynomial of degree $n$, then
\begin{align}
\left\{\dfrac{1}{2\pi}\int_0^{2\pi}|P'(e^{i\theta})|^pd\theta\right\}^\frac{1}{p}\le
n\left\{\dfrac{1}{2\pi}\int_0^{2\pi}|P(e^{i\theta})|^pd\theta\right\}^\frac{1}{p}
\end{align}
Equality holds in (1.2) for $P(z)=\alpha z^n$, $|\alpha|\neq 0$. If
we let $p\to \infty$ in (1.2), we get inequlaity (1.1).

Let $\alpha$ be a complex number. If $P(z)$ is a polynomial of degree $n$, then the polar derivative of $P(z)$
with respect to the point $\alpha$, denoted by $D_\alpha P(z)$, is defined by
\[
D_\alpha P(z)=nP(z)+(\alpha -z)P'(z)
\]
clearly $D_\alpha P(z)$ is a polynomial of degree at most $n-1$ and it generalizes the ordinary derivative in the sense that
\begin{align}
\lim_{\alpha\to\infty}\left[\dfrac{D_\alpha P(z)}{\alpha}\right]=P'(z).
\end{align}
As an extension of (1.1) to the polar derivative, Aziz and Shah \cite{3}, have shown that if $P(z)$ is a polynomial of degree $n$,
then for every real or complex number $\alpha$ with $|\alpha|>1$ and for $|z|=1$,
\begin{align}
|D_\alpha P(z)|\le n|\alpha|\max\limits_{|z|=1} |P(z)|
\end{align}
As a generalization of (1.2) to the polar derivaative Aziz et al.\ \cite{4}, proved the following result.

\begin{theoremA}\rm
If $P(z)$ is a polynomial of degree $n$, then for every complex number $\alpha$ with $|\alpha|\ge 1$ and $p\ge 1$
\begin{align}
\left\{\int_0^{2\pi}|D_\alpha P(e^{i\theta})|^p d\theta\right\}^{\frac{1}{p}}\le n(|\alpha|+1)
\left\{\int_0^{2\pi}|P(e^{i\theta})|^p d\theta\right\}^{\frac{1}{p}}
\end{align}
For the Class of polynomials having no zeros in $|z|<1$, inequality (1.2) was improved by D-Bruijin \cite{5} that if $P(z)\neq 0$ in $|z|<1$,
then for $p\ge 1$
\begin{align}
\left\{\int_0^{2\pi}|P'(e^{i\theta})|^p d\theta\right\}^{\frac{1}{p}}\le n C_p
\left\{\int_0^{2\pi}|P(e^{i\theta})|^p d\theta\right\}^{\frac{1}{p}}
\end{align}
where
\begin{align}
c_p=\left\{\dfrac{1}{2\pi}\int_0^{2\pi}|1+e^{i\theta}|^p d\theta\right\}^{-\frac{1}{p}}
\end{align}
\end{theoremA}

As an extension to the polar derivative. A. Aziz and N. Rather \cite{6}, proved the following generalization of (1.5). In fact they proved.

\begin{theoremA}\rm
If $P(z)$ is a polynomial of degree $n$ which does not vanish in $|z|<1$, then for every complex number $\alpha$ with
$|\alpha|\ge 1$ and $p\ge 1$
\begin{align}
\left\{\int_0^{2\pi}|D_\alpha P(e^{i\theta})|^p d\theta\right\}^{\frac{1}{p}}\le n (|\alpha|+1)c_p
\left\{\int_0^{2\pi}|D_\alpha P(e^{i\theta})|^p d\theta\right\}^{\frac{1}{p}}
\end{align}
where $C_p$ is defined by (1.7).
\end{theoremA}

Recently, Irshad et al.\ \cite{7} proved the following result.

\begin{theoremA}\rm
If $P(z)$ is a polynomial of degre $n$ which does not vanish in $|z|<K\le 1$, then for every
$\alpha,\beta\in C$ with $|\alpha|\ge K$, $|\beta|\le 1$ and $p\ge 1$
\begin{align}
&\left\{\int_0^{2\pi}\left|e^{i\theta}D_\alpha P(e^{i\theta})+n\dfrac{(|\alpha|-K)}{K+1} \beta P(e^{i\theta})\right|^pd\theta\right\}^{\frac{1}{p}}\nonumber\\
&\quad \le  n \left(1+|\alpha|+2\dfrac{(|\alpha|-K)}{K+1}|\beta|\right)C_p \left\{\int_0^{2\pi}|P(e^{i\theta})|^p d\theta\right\}^{\frac{1}{p}}
\end{align}
where $C_p$ is defined by (1.7).
\end{theoremA}

In this paper we prove the following more general result which also generalize Theorem~B
and yields a number of known polynomial inequalities.

\def\eit{e^{i\theta}}
\def\eip{e^{i\phi}}
\begin{theorem}\rm
If $P(z)=a_nz^n+\sum\limits_{j=\mu}^n a_{n-j}z^{n-j}$, $1\le\mu\le n$ be a polynomial of degree $n$
which does not vanish in $|z|<K\le 1$, then for every $\alpha,\beta\in C$ with $|\alpha|\ge K$, $|\beta|\le 1$ and $p\ge 1$
\begin{align}
&\left\{\int_0^{2\pi}\left|e^{i\theta}D_\alpha P(e^{i\theta})+n\dfrac{(|\alpha|-K^\mu)}{K^\mu+1} \beta P(e^{i\theta})\right|^pd\theta\right\}^{\frac{1}{p}}\nonumber\\
&\quad \le  n \left(1+|\alpha|+2\dfrac{(|\alpha|-K^\mu)}{K^\mu+1}|\beta|\right)C_p \left\{\int_0^{2\pi}|P(e^{i\theta})|^p d\theta\right\}^{\frac{1}{p}}
\end{align}
where $C_p$ is defined by (1.7), or equivalently
\begin{align}
&\left\|e^{i\theta}D_\alpha P(\eit)+n\dfrac{(|\alpha-K^\mu|)}{K^\mu+1}\beta P(\eit)\right\|_p\nonumber\\
&\quad \le n\left(1+|\alpha|+2\dfrac{(|\alpha|-K^\mu)}{K^\mu+1}|\beta|\right)\dfrac{\|P(\eit)\|_p}{\|1+\eip\|_p}
\end{align}
\end{theorem}

\begin{remark}
If we choose $\mu=1$ in (1.10), we get Theorem~C and if we choose $\beta=0$ and $K=1$ in (1.10), we get Theorem~B.
\end{remark}

\section{Lemmas}

For the proof of this theorem, we need the following lemmas. The first lemma is due to Gulshan  Singh et al. \cite{8}.

\begin{lemma}\rm
Let $P(z)=a_n z^n+\sum\limits_{j=\mu}^n a_{n-j}z^{n-j}$, $1\le \mu\le n$ be a polynomial of degree having all its zeros
in the disk $|z|\le K$, $K\le 1$, then for every real
or complex number $\alpha$ with $|\alpha|\ge K$, $K\le 1$ and for $|z|=1$
\[
|D_\alpha P(z)|\ge n\left(\dfrac{|\alpha|-K^\mu}{K^\mu}\right)|P(z)|
\]
\end{lemma}

\begin{lemma}\rm
Let $Q(z)$ be a polynomial of degree $n$ having all its zeros in $|z|<K$, $K\le 1$ and $P(z)$ is a polynomial of degree not exceeding that of $Q(z)$.
If $|P(z)|\le |Q(z)|$ for $|z|=K\le 1$, then for every $\alpha,\beta \in C$ with $|\alpha|\ge  K$, $\beta|\le 1$
\[
\left|zD_\alpha P(z)+n\dfrac{(|\alpha|-K^\mu)}{K^\mu+1}\beta P(z)\right|
\le \left|zD_\alpha Q(z)+n\dfrac{(|\alpha|-K^\mu)}{K^\mu+1}\beta Q(z)\right|
\]
\end{lemma}

\proof
Since $|\lambda P(z)|\le |P(z)|\le |Q(z)|$, for $\lambda <1$ and $|z|=K$, then for Rouche's Theorem $Q(z)-\lambda P(z)$ and $Q(z)$ have the same number of
zeros in $|z|<K$.
On the other hand by inequality $|P(z)|\le |Q(z)|$ for $|z|=K$, any zero of $Q(z)$, that lies on $|z|=K$, in the zero of $P(z)$. Therefore, $Q(z)-\lambda P(z)$ has all its zero in the closed disk $|z|\le K$. Hence by Lemma~1 for all real or complex numbers $\alpha$ with $|\alpha|\ge K$ and $|z|=1$, we have
\begin{align}
|zD_\alpha(Q(z)-\lambda P(z))|\ge n\dfrac{(|\alpha|-K^\mu)}{K^\mu+1}|Q(z)-\lambda P(z)|
\end{align}
Now consider a similar argument that for any value of $\beta$ with $|\beta|<1$, we have
\begin{align}
|zD_\alpha(Q(z)-\lambda P(z))|
&\ge n\left(\dfrac{|\alpha|-K^\mu}{K^\mu+1}\right)|Q(z)-\lambda P(z)|\nonumber\\
&>   n|\beta|\left(\dfrac{|\alpha|-K^\mu}{K^\mu+1}\right)|Q(z)-\lambda P(z)|
\end{align}
where $|z|=1$, resulting in
\begin{align}
T(z)=|zD_\alpha Q(z)-\lambda zD_\alpha P(z)|+n\beta
\dfrac{|\alpha|-K^\mu}{K^\mu+1}\{Q(z)-\lambda P(z)\}\neq 0
\end{align}
where $|z|=1$.

That is
\begin{align}
T(z)=\left|zD_\alpha Q(z)+n\beta\dfrac{|\alpha|-K^\mu}{K^\mu+1}Q(z)\right|
-\lambda \left|zD_\alpha P(z)+n\beta\dfrac{|\alpha|-K^\mu}{K^\mu+1}P(z)\right|\neq 0
\end{align}
for $|z|=1$

We also conclude that
\begin{align}
\left|zD_\alpha Q(z)+n\beta\dfrac{|\alpha|-K^\mu}{K^\mu+1}Q(z)\right|
\ge \left|zD_\alpha P(z)+n\beta\dfrac{|\alpha|-K^\mu}{K^\mu+1}P(z)\right|
\end{align}
for $|z|=1$.

If (2.5) is not true, then there is a point $z=z_0$ with $|z_0|=1$, such that
\begin{align}
\left|z_0D_\alpha Q(z_0)+n\beta\dfrac{|\alpha|-K^\mu}{K^\mu+1}Q(z_0)\right|
< \left|z_0D_\alpha P(z_0)+n\beta\dfrac{|\alpha|-K^\mu}{K^\mu+1}P(z_0)\right|
\end{align}
Take
\begin{align}
\lambda=\dfrac{z_0D_\alpha Q(z_0)+n\beta\dfrac{|\alpha|-K^\mu}{K^\mu+1}Q(z_0)}
{z_0D_\alpha P(z_0)+n\beta\dfrac{|\alpha|-K^\mu}{K^\mu+1}P(z_0)}
\end{align}
then $|\lambda|<1$ with this choice, we have from (2.4),
$T(z_0)=0$ for $|z_0|=1$. But this contradicts the fact that $T(z)\neq 0$ for $|z|=1$.\
For $\beta$ with $|\beta|=1$, (2.5) follows by continuity.

This completes the proof.

The next lemma is due to Aziz and Rather \cite{4}.

\begin{lemma}\rm
If $P(z)$ is a polynomial of degree $n$ such that $P(0)\neq 0$ and $Q(z)=z^n \overline{p\left(\dfrac{1}{\bar{z}}\right)}$, then for every $p\ge 0$ and $\phi$ real
\[
\int_0^{2\pi}\int_0^{2\pi} |Q'(\eit)+\eit P'(\eit)|^p d\theta d\phi\le n^p \int_0^{2\pi} |P(\eit)|^p d\theta
\]
\end{lemma}

\section{Proof of the theorem}

\noindent\textbf{Proof of the Theorem.}
Let $P(z)$ be a polynomial of degree $n$ which does not vanish in $|z|<K$, $K\le 1$.
By Lemma~2 for complex numbers $\alpha$, $\beta$ with $|\alpha|\ge K$, $|\beta|\le 1$, we have
\begin{align}
\left|zD_\alpha P(z)+n\dfrac{(|\alpha|-K^\mu)}{K^\mu+1}\beta P(z)\right|
\le \left|zD_\alpha Q(z)+n\dfrac{(|\alpha|-K^\mu)}{K^\mu+1}\beta Q(z)\right|
\end{align}
For every real $\phi$ and $\xi\ge 1$, we have
\[
|\xi+\eit|\ge |1+\eip|
\]
which implies for any $p\ge 0$
\begin{align}
\left\{\int_0^{2\pi}|\xi+\eip|^p d\phi\right\}^\frac{1}{p} \ge \left\{\int_0^{2\pi}|1+\eip|^p d\phi \right\}^\frac{1}{p}
\end{align}
If $\eit D_\alpha P(\eit)+n\dfrac{(|\alpha|-K^\mu)}{K^\mu+1}\beta P(\eit)\neq 0$, we can take
\[
\xi=\dfrac{\eit D_\alpha Q(\eit)+n \dfrac{(|\alpha|-K^\mu)}{K^\mu+1}\beta Q(\eit)}{\eit D_\alpha P(\eit)+n \dfrac{(|\alpha|-K^\mu)}{K^\mu+1}\beta P(\eit)}
\]
where according to (3.1), $|\xi|\ge 1$. Now
\begin{align}
&\int_0^{2\pi}\left|\eit D_\alpha Q(\eit)+n\dfrac{(|\alpha|-K^\mu)}{K^\mu+1} \beta Q(\eit)\right.\nonumber\\
&\qquad \left.+\eip\left[\eit D_\alpha P(\eit)+n\dfrac{(|\alpha|-K^\mu)}{K^\mu+1} \beta P(\eit) \right]\right|^p d\phi\nonumber\\
&=\left|\eit D_\alpha P(\eit)+n \dfrac{(|\alpha|-K^\mu)}{K^\mu+1} \beta P(\eit)\right|^p \int_0^{2\pi} |\xi+\eip|^p d\phi\nonumber\\
&\ge \left|\eit D_\alpha P(\eit)+n \dfrac{(|\alpha|-K^\mu)}{K^\mu+1} \beta P(\eit)\right|^p \int_0^{2\pi} |1+\eip|^p d\phi
\end{align}
This inequality is trivially true if
\[
\eit D_\alpha P(\eit)+n \dfrac{(|\alpha|-K^\mu)}{K^\mu+1} \beta P(\eit)=0
\]
Integrating both sides of (3.3) with respect to $\theta$ from 0 to $2\pi$, we have
\begin{align}
&\int_0^{2\pi}\int_0^{2\pi}\left|\eit D_\alpha Q(\eit)+n\dfrac{(|\alpha|-K^\mu)}{K^\mu+1} \beta Q(\eit)\right.\nonumber\\
&\qquad \left.+\eip\left[\eit D_\alpha P(\eit)+n\dfrac{(|\alpha|-K^\mu)}{K^\mu+1} \beta P(\eit) \right]\right|^p d\theta\,d\phi\nonumber\\
&\ge \int_0^{2\pi} \left|\eit D_\alpha P(\eit)+n \dfrac{(|\alpha|-K^\mu)}{K^\mu+1} \beta P(\eit)\right|^p d\theta \int_0^{2\pi} |1+\eip|^p d\phi
\end{align}
Now for $0\le \theta <2\pi$
\begin{align}
&\left|\eit D_\alpha Q(\eit)+n\dfrac{(|\alpha|-K^\mu)}{K^\mu+1} \beta Q(\eit) +\eip \left[\eit D_\alpha P(\eit)+n \dfrac{(|\alpha|-K^\mu)}{K^\mu+1} \beta P(\eit)\right]\right|\nonumber\\
&\quad=\left|\left[\eit\{nQ(\eit)+(\alpha-\eit)Q'(\eit)\}+n\dfrac{(|\alpha|-K^\mu)}{K^\mu+1} \beta Q(\eit)\right]\right.\nonumber\\
&\qquad+\left.\eip\left[\eit\{nP(\eit)+(\alpha-\eit)P'(\eit)\}+n\dfrac{(|\alpha|-K^\mu)}{K^\mu+1} \beta P(\eit)\right]\right|
\end{align}
\begin{align}
&\quad=\left|\left[\eit\{nQ(\eit)-\eit Q'(\eit)\}+\alpha\eit Q'(\eit)+n\dfrac{(|\alpha|-K^\mu)}{K^\mu+1} \beta Q(\eit)\right]\right.\nonumber\\
&\qquad+\left.\eip\left[\eit\{nP(\eit)-\eip P'(\eit)\}+\alpha \eit P'(\eit)+n\dfrac{(|\alpha|-K^\mu)}{K^\mu+1} \beta P(\eit)\right]\right|
\end{align}
Since $Q(z)=z^n\overline{P\left(\dfrac{1}{\bar{z}}\right)}$, we have $P(z)=z^n \overline{Q\left(\dfrac{1}{\bar{z}}\right)}$ and it can be easily verified that for $0\le \theta <2\pi$
\[
nP(\eit)-\eit P'(\eit)=e^{i(n-1)\theta}\overline{Q'(\eit)}
\]
and
\[
nQ(\eit)-\eit Q'(\eit)=e^{i(n-1)\theta}\overline{P'(\eit)}
\]
From (3.6)
\begin{align}
&\left|\eit D_\alpha Q(\eit)+n\dfrac{(|\alpha|-K^\mu)}{K^\mu+1} \beta Q(\eit) +\eip \left[\eit D_\alpha P(\eit)+n \dfrac{(|\alpha|-K^\mu)}{K^\mu+1} \beta P(\eit)\right]\right|\nonumber\\
&\quad=\left|\left[\eit\{ e^{i(n-1)\theta} \overline{P'(\eit)}\}\right]+\alpha \eit \left[Q'(\eit)+\eip P'(\eit)\right]\right.\nonumber\\
&\qquad\left.+n\dfrac{(|\alpha|-K^\mu)}{K^\mu+1}\beta [Q(\eit)+\eip P(\eit)]+\eip \eit e^{i(n-1)\theta}Q'(\eit)\right|
\end{align}
Therefore, (3.4) in conjunction with (3.7) gives
\begin{align}
&\left\{\int_0^{2\pi}\int_0^{2\pi} |\eit e^{i(n-1)\theta}\{\overline{P'(\eit)}+\eip \overline{Q'(\eit)}\} +\alpha \eip [Q'(\eit)\eip P'(\eit)]\right.\nonumber\\
&\quad\left.+n\dfrac{(|\alpha|-K^\mu)}{K^\mu+1}\beta [Q(\eit)+\eip P(\eit)]|^p d\theta d\phi\right\}^{\frac{1}{p}}\nonumber\\
&\quad\ge \left\{\int_0^{2\pi} \left|\eit D_\alpha P(\eit)+n\dfrac{(|\alpha|-K^\mu)}{K^\mu+1}\beta P(\eit) \right|^p d\theta \int_0^{2\pi} |1+\eip|^p d\phi\right\}^\frac{1}{p}
\end{align}
By Minkowski inequality, we have
\begin{align*}
&\left\{\int_0^{2\pi} |\eit D_\alpha P(\eit)+n\dfrac{(|\alpha|-K^\mu)}{K^\mu+1} \beta P(\eit)|^p d\theta \int_0^{2\pi} |1+\eip|^p d\phi \right\}^\frac{1}{p}\nonumber\\
&\quad\ge \left\{\int_0^{2\pi}\int_0^{2\pi} |Q'(\eit)+\eip P'(\eit)|^p d\theta d\phi \right\}^{\frac{1}{p}}\{1+|\alpha|\}\nonumber\\
&\qquad + \left|n\dfrac{(|\alpha|-K^\mu)}{K^\mu+1}\beta\right|\left\{\int_0^{2\pi}\int_0^{2\pi} |Q(\eit)+\eip P(\eit)|^p d\theta d\phi \right\}^{\frac{1}{p}}
\end{align*}
By Lemma~3, we have
\begin{align*}
&\left\{\int_0^{2\pi} |\eit D_\alpha P(\eit)+n\dfrac{(|\alpha|-K^\mu)}{K^\mu+1} \beta P(\eit)|^p d\theta \int_0^{2\pi}|1+\eip|^p d\phi \right\}^{\frac{1}{p}}\nonumber\\
&\quad\le \left\{ 2n^p\pi \int_0^{2\pi}|P(\eit)|^p d\theta\right\}^{\frac{1}{p}}\{1+|\alpha|\}\nonumber\\
&\qquad + 2n\left|\dfrac{(|\alpha|-K^\mu)}{K^\mu+1}\right|\left\{2\pi\int_0^{2\pi}|P(\eit)|^p d\theta \right\}^{\frac{1}{p}}\\
&=\left[n(1+|\alpha|)+2n\dfrac{(|\alpha|-K^\mu)}{K^\mu+1}|\beta|\right]\left\{2\pi \int_0^{2\pi} |P(\eit)|^p d\theta\right\}^\frac{1}{p}
\end{align*}
This implies
\begin{align*}
&\left\{\int_0^{2\pi} |\eit D_\alpha P (\eit)+n\dfrac{(|\alpha|-K^\mu)}{K^\mu+1} \beta P(\eit)|^p d\theta \right\}^{\frac{1}{p}}\nonumber\\
&\le n\left(1+|\alpha|+2\dfrac{(|\alpha|-K^\mu)}{K^\mu+1}|\beta|\right)C_p\left\{\int_0^{2\pi} |P(\eit)|^p d\theta\right\}^\frac{1}{p}
\end{align*}
where $C_p$ in defined by (1.7),

\noindent or equivalently,
\begin{align*}
&\left\|e^{i\theta}D_\alpha P(\eit)+n\dfrac{(|\alpha|-K^\mu)}{K^\mu+1}\beta P(\eit)\right\|\nonumber\\
&\quad \le n\left(1+|\alpha|+2\dfrac{(|\alpha|-K^\mu)}{K^\mu+1}|\beta|\right)\dfrac{\|P(\eit)\|_p}{\|1+P(\eip)\|_p}
\end{align*}
Hence the result.

\end{document}